\documentclass[11pt]{article}
\usepackage{amsmath}
\usepackage{graphicx}
\usepackage{amsfonts, amsmath, amsthm, amssymb, latexsym, amscd, eufrak}

\newtheorem{thm}{\bf Theorem}
\newtheorem{lem}{\bf Lemma}

\newtheorem{exam}{\bf Example}

\newenvironment{prf}{{\noindent\bf Proof}\rm}{\hfill{$\Box$}}

\begin{document}

\title{The balanced 2-median and 2-maxian problems on a tree}


\author{
\bf{Jafar Fathali}\thanks{Corresponding author} \\
      Faculty of Mathematical Science, Shahrood University of Technology,\\
      University Blvd., Shahrood, Iran, email: fathali@shahroodut.ac.ir\\ \\
\bf{Mehdi Zaferanieh}\\
       Department of Mathematics, Hakim Sabzevari University, Tovhid town,\\
       Sabzevar, Iran, email: mehdi.zaferanieh@gmail.com
 }

\date{}
\maketitle
\newpage

\begin{abstract}
This paper deals with the facility location problems with balancing on allocation clients to servers. Two bi-objective models are considered, in which one objective is the traditional p-median or p-maxian objective and the second is to minimize the maximum demand volume allocated to any facility.
An edge deletion method with time complexity $O(n^2)$ is presented for the balanced $2$-median problem on a tree. For the balanced $2$-maxian problem, it is shown the optimal solution is two end vertices of the diameter of the tree, which can be obtained in a linear time.
\end{abstract}

\begin{quote}

\end{quote}

\section{Introduction}

The most location problems are concerned about minimization of transportation time between servers and clients. However in real applications usually the service time is as considerable as transportation time due to attracting clients. So in this paper we consider a bi-objective problem corresponding to the transportation time as well as service time. The $p$-median and $p$-maxian objective functions have been considered for transportation time while an avoidance congestion function has been considered for service time. We call these problems balanced $p$-median and $p$-maxian problems, respectively. To balance the facility servicing times, we propose minimizing the maximum number of clients that are served by facilities.

Let $G=(V,E)$ be a given graph, where $V$ is the set of vertices and $E$ is the set of edges. Let $|V|=n$ and $|E|=m$. The $p$-median problem asks to find a set of $p$ vertices of $G$, called facilities, such that the sum of weighted distances from vertices to the closest facility is minimized. Kariv and Hakimi \cite{KH79} showed that the $p$-median problem is $NP$-hard on general networks while it can be solved in polynomial time on tree networks. The initial works of the $p$-median problem is referred to Hakimi \cite{H64,H65}. Hakimi \cite{H65} showed that at least one optimal solution of the $p$-median problem is located on vertices. Kariv and Hakimi \cite{KH79} presented an $O(p^2n^2)$ time algorithm for this problem on a tree. The time complexity is improved to $O(pn^2)$ by Tamir \cite{T96}. In the case $p=2$ on a tree, Gavish and Sridhar \cite{GS95} presented an $O(nlogn)$ algorithm.

The obnoxious case of the $p$-median problem is called $p$-maxian problem. In the $p$-maxian problem a set which contains $p$ vertices is sought so that the sum of weighted distances of clients to the farthest facility is maximized. The NP-hardness of this problem on general networks is shown in \cite{HM88}. Zelinka~\cite{Z68} showed that in the tree graphs an optimal solution of the $1$-maxian problem is contained on the leaves. Ting~\cite{T84} proposed a linear time algorithm to the $1$-maxian problem. The $1$-maxian problem on general networks is investigated by Church and Garfinkel \cite{CG78}. They presented an $O(mn\log n)$ time algorithm, and  Tamir \cite{T91} improved the time complexity to $O(mn)$.  Burkard et al. \cite{BFK07} showed that the optimal solution of the $2$-maxian problem on a tree lies on the two end vertices of the diameter, where diameter is the longest path of the tree. They also showed that $p$-maxian problem on the tree is reduced to the $2$-maxian. Based on these properties they presented a linear time algorithm for the $p$-maxian problem on a tree. Kang and Cheng \cite{KC10} extended the algorithm to the case that the underlying network is a block graph.

The equity location models are introduced in the last two decades. In these models, facilities are to be located to maximize the equality between the demand points. Some researchers have been attracted to this subject. Among them Gavalec and Hudec \cite{GH95} considered an equity model which its objective function is the maximum difference in the distance from a demand point to its farthest and nearest facility. They called this problem as balancing function. Berman et al. \cite{BDTW09} considered the problem of finding the location of $p$ facilities such that the weights attracted to the different facilities are as close as possible. They formulated this problem as minimizing the maximum weight assigned to each facility.  Marin \cite{M11} considered the balanced discrete location problem, in which the objective function is minimizing the difference between the maximum and the minimum weights allocated to different facilities. Barbati and Piccolo \cite{BP15} proposed some properties to describe the behavior of the equality measures in an optimization context.

Another paper related to balanced facility location models is that of Lejeune and Prasad \cite{LP13}, which propose models to investigate effectiveness-equity tradeoffs in tree network facility location problems. The 1-median objective is considerd as a measure of effectiveness, and the Gini index is used as a measure of equity. A bi-criteria problems involving these objectives is presented. Landete and Marin \cite{LM14}, considered the spanning trees with balanced weights, i.e., where the differences among the weights are minimized.

In the paper of López-de-los-Mozos et al. \cite{LMP08} the ordered weighted averaging operator is applied to define a model which generalizes some inequality measures. In their work,  for a location $x$, the value of the objective function is the ordered weighted average of the absolute deviations from the average distance from the facilities to the location $x$.
We refer the interested reader to \cite{MS94,EL95}, two reviews of the literature on equity measurement in location theory.

In the next section, we formulate the balanced $p$-median and $p$-maxian problems. In Section 3, an $O(n^2)$ time algorithm for balanced 2-median problem on a tree is presented. The balanced $2$-maxian problem on a tree is investigated in Section 4, and a linear time algorithm is proposed for this problem.

\section{Problem definition}
Let $G=(V,E)$ be a graph where $V$ is the set of vertices and $E$ is the set of edges and $|V|=n$. Each vertex $v_i$ has a non-negative weight $w_i$, which is the number of clients on vertex $v_i$. The weight $w_i$ also called the demand at $v_i$. Let $d_{ij}$ be the distance between vertices $v_i$ and $v_j$. The $p$-median problem asks to find a set of $p$ vertices, $X_p=\{x_1,...,x_p\}$ called facilities, so that the sum of the weighted distances from the vertices to the closest facility in $X_p$ is minimized, i.e.
$$\min f_1(X_p)=\sum_{i=1}^{n}w_id(X_p,v_i),$$
where for any vertex $v\in V$, $d(X_p,v)=\min_{x_j\in X_p}d(x_j,v).$

In the $p$-maxian problem the goal is maximizing the sum of the weighted distances from the vertices to the farthest facility in $X_p$, i.e.
$$\max f_2(X_p)=\sum_{i=1}^{n}w_i\max_{x_j\in X_p} d(x_j,v_i).$$

We define the bi-objective models of balanced $p$-median and $p$-maxian problems as follows. Let $t_i$ be the service time by different facilities for serving clients on vertex $v_i$. Also let $V_j$ be the set of vertices of $V$ that are allocated to the facility $x_j$. Then, in the balanced $p$-median problem we should find a partition $\{G_1=(V_1,E_1),...,G_p=(V_p,E_p)\}$ of $G$ such that,

\begin{align}
&\min f_1(X_p)\\
&\min f_3(X_p)=\max\{\sum_{v_i\in V_j}w_i t_i,~~j=1,...,p\} \label{ps}.
\end{align}


Similarly, in the bi-objective function of balanced $p$-maxian problem we would find a partition of $G$ such that,

\begin{align}
&\max f_2(X_p)\\
&\min f_3(X_p)=\max\{\sum_{v_i\in V_j}w_i t_i,~~j=1,...,p\} \nonumber.
 \end{align}

Note that the goal in these problems is partitioning graph $G=(V,E)$ into $p$ subgraphs. Then in the balanced $p$-median problem, in each partition the median should be determined. While in the balanced $p$-maxian problem, for $i=1,...,p$, each vertex $u\in V_i$ is allocated to the facility $x_i \in X_p $ that $d(u,x_i)=\max\{d(u,x_j)|x_j\not\in V_i,~j=1,...,p\}$.
Therefor, by this partitioning, we can write,

$$f_1(X_P)=f_2(X_p)=\sum_{j=1}^{p}\sum_{v_i\in V_j}w_id(x_j,v_i).$$

In the classical $p$-median problem each client is allocated to the closest facility and $x_i\in V_i$ for $i=1,...,p$, while in the $p$-maxian problem each vertex is allocated to the farthest facility and $x_i\not\in V_i$ for $i=1,...,p$. However, in the balanced case either client may be allocated to the closest or farthest facility. 

Let $z_i=w_it_i$ for $i=1,...,n$ then
\begin{align}
f_3(X_p)=\max\{\sum_{v_i\in V_j}z_i,~~j=1,...,p\}.
\end{align}

In the case that the clients' service times are all equal, i.e. $t_i=t\geq 0$ for $i=1,...,n$, then the objective function in (\ref{ps}) reduces to the following:
\begin{align}
\min f_{4}(X_p)=\max\{\sum_{v_i\in V_j}w_i,~~j=1,...,p\}.
\end{align}

The problem of minimizing $f_4(X)$ is considered by Berman et al. \cite{BDTW09}. They showed this problem is $NP$-hard. Since the $p$-median and $p$-maxian problems are also $NP$-hard, then the balanced cases are $NP$-hard, too.
%
%
%
%
%

%

In this paper, we use the weighted sum method to the proposed bi-objective problems.
%
The weighted sum problem of bi-objective $p$-median and $p$-maxian problems, can be written as
$$\min f_{pmed}(X)=\lambda f_1(X)+(1-\lambda)f_3(X),$$
$$\max f_{pmax}(X)=\lambda f_2(X)-(1-\lambda)f_3(X),$$
where $0\leq\lambda\leq 1$. Note that in these models $1-\lambda$ can be interpreted as the servers' balanced coefficient.

Next, the balanced $2$-median and $2$-maxian problems are studied.

\begin{lem}\label{f2}
In the case $p=2$, the model $f_3(X_p)$ can be represented as the following problem,
\begin{align}f_5(X_p)=|\sum_{v_i\in V_1}z_i -\sum_{v_i\in V_2}z_i|.\end{align}
\end{lem}

\begin{prf}
Let $Z=\sum_{i=1}^nz_i$, then in the case $p=2$, the objective function $f_3(X_p)=\max\{\sum_{v_i\in V_j}z_i,~~j=1,...,p\}$ can
be represented as
$$f_3(X_p)=\max\{\sum_{v_i\in V_1}z_i,~~\sum_{v_i\in V_2}z_i\}$$
$$=\frac{|\sum_{v_i\in V_1}z_i +\sum_{v_i\in V_2}z_i|+|\sum_{v_i\in V_1}z_i -\sum_{v_i\in V_2}z_i|}{2}$$
$$=\frac{Z+|\sum_{v_i\in V_1}z_i -\sum_{v_i\in V_2}z_i|}{2}.$$
Since $Z$ is a fixed value, it can be left out from the objective function.
\end{prf}

Obviously, by Lemma \ref{f2} in the case $p=2$, any partition of $V$ to $V_1$ and $V_2$ such that $\sum_{v_i\in V_1}z_i=\sum_{v_i\in V_2}z_i$ provides the optimal solution of $\min f_3(X_p)$.


\section{Balanced 2-median problem on a tree}

In this section, the balanced 2-median problem on a tree is considered. We propose the following weighted objective function to find the Pareto optimal solutions to the balanced 2-median problem.

\begin{align}
\min f_{pmed}(X)=\lambda(\sum_{v_i\in T_1}w_i d(v_i,x_1)+\sum_{v_i\in T_2}w_i d(v_i,x_2))+(1-\lambda)|\sum_{v_i\in T_1}z_i-\sum_{v_i\in T_2}z_i|. \label{bpmed}
\end{align}
Where $T$ is the underlaying network and $T_1$ and $T_2$ are the two subtrees of $T$ which contain vertices in $V_1$ and $V_2$, respectively.


A well-known method for solving the classical 2-median problem on a tree is the edge deletion method (see e.g.\cite{BCD00}). In this method by deleting any edge the median of obtained subtrees are found and the best one is chosen as the solution of the 2-median problem. Since in the second part of our model partitioning the tree in two subtrees is necessary, we apply the edge deletion method to solve the balanced 2-median problem as given below.

For every edge $e$, let $T_1^e$ and $T_2^e$ be two subtrees of $T$ which are obtained by deleting edge $e$. Let $m_1^e$ and $m_2^e$ be the 1-medians of $T_1^e$ and $T_2^e$, respectively. We calculate the value of objective function $f_{pmed}(.)$ as follows:
\begin{align}
f_{pmed}(m_1^e,m_2^e)=\lambda(\sum_{v_i\in T_1^e}w_i d(v_i,m_1^e)+\sum_{v_i\in T_2^e}w_i d(v_i,m_2^e)) +(1-\lambda)|\sum_{v_i\in T_1^e}z_i -\sum_{v_i\in T_2^e}z_i|.\label{btr}
\end{align}
Then the pairs of medians corresponding to the minimum amounts of Problem (\ref{btr}) are chosen as the solution of the balanced 2-median problem. The time complexity of this method is $O(n^2)$.

In the following example we show that in the optimal solution of Problem (\ref{bpmed}) the customers may not be assigned to the closest median.

\begin{exam}\label{nonear}
Consider the tree depicted in Fig. \ref{fig2}, where the weights and service times of all vertices are equal to one. The solution of the balanced 2-median problem for the case $\lambda=\frac{1}{2}$ is $\{v_2,v_4\}$. If we allocate $v_3$ to $v_4$ then the value of objective function is $\frac{1}{2}(4+2)$ while allocating $v_3$ to $v_2$ results in  the value of objective function equal to $\frac{1}{2}(5+0)$. Therefore, the optimal solution is obtained by deleting edge $(v_3,v_4)$.

\begin{figure}
\vspace{1.2cm}
\centerline{\includegraphics[width=10cm]{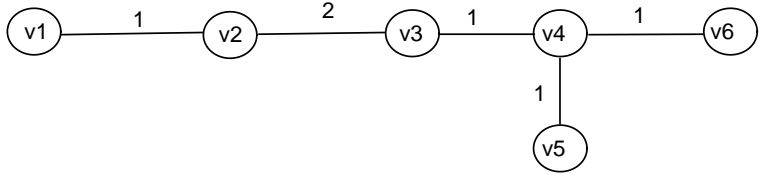}}
\caption{\label{fig2}\small The tree for Example \ref{nonear}.}
\end{figure}
\end{exam}

\section{Balanced 2-maxian problem on a tree}
In this section, we consider the balanced 2-maxian problem on the tree $T$. Let $T_1$ and $T_2$ be the two subtrees of $T$ which contain vertices in $V_1$ and $V_2$, respectively. The model of this problem is given as follows:
\begin{align}\max f_{pmax}(X)=\lambda(\sum_{v_i\in T_1}w_i d(v_i,x_1)+\sum_{v_i\in T_2}w_i d(v_i,x_2))-(1-\lambda)|\sum_{v_i\in T_1}z_i-\sum_{v_i\in T_2}z_i|.\end{align}

Burkard et al. \cite{BFK07} showed that the solution of the $2$-maxian problem is two end vertices of a diameter which can be obtained in a linear time. For $\lambda\in[0,1]$ the balanced 2-maxian problem also could be solved by such edge deletion method as in Burkard et al. \cite{BCD00} which is presented for the $2$-median problem on a tree with positive and negative weights.

\subsection{The edge deletion method}

Let $T_a$ and $T_b$ be two subtrees of $T$ which are obtained by deleting edge $e=(a,b)$. We construct two trees $T_1^e$ and $T_2^e$ which their vertices are the same as $T$. The weight of vertices in $T_1^e$ are defined as
$$ w^1_i=\left\{
  \begin{array}{l l}
w_i &  \text{if}~ v_i \in T_b \\
0 & otherwise,
 \end{array}\right.$$
and analogously the weights of vertices in $T_2^e$ are defined as
$$ w^2_i=\left\{
  \begin{array}{l l}
w_i &  \text{if}~ v_i \in T_a \\
0 & otherwise.
 \end{array}\right.$$
Then the set $X=\{x_1,x_2\}$ is considered which maximizes the following objective function.
$$f_{pmax}(X)=\lambda(\sum_{v_i\in T^e_1}w_i d(v_i,x_1)+\sum_{v_i\in T^e_2}w_i d(v_i,x_2))-(1-\lambda)|\sum_{v_i\in T^e_1}z_i-\sum_{v_i\in T^e_2}z_i|.$$
Note that for each edge we should compute the objective function for any pair $\{x_1,x_2\}$. Therefore, the time complexity of this method is $O(n^3)$.

The following example shows in the optimal solution of the balanced $2$-maxian problem, clients may not be allocated to the farthest facility.

\begin{exam}\label{nofar}
Consider again the tree in Fig. \ref{fig2}, where the weights of edges $(v_1,v_2)$ and $(v_2,v_3)$ are replaced by 2 and 1, respectively. The weights and service times of all vertices are equal to one. The solution of balanced 2-maxian problem for the case $\lambda=\frac{1}{2}$ is $\{v_1,v_6\}$ (and $\{v_1,v_5\}$) which is obtained by deleting edge $(v_3,v_4)$. The vertices $v_4$, $v_5$ and $v_6$ are allocated to the facility in $v_1$, and the vertices $v_1$, $v_2$ and $v_3$ are allocated to the facility in $v_6$. Note that the vertex $v_3$ is allocated to facility in $v_6$, but the farthest facility from $v_3$ is in $v_1$. The optimal value of objective function is $\frac{1}{2}(24-0)$. If the vertex $v_3$ is allocated to the facility in $v_1$ then the value of objective function is  $\frac{1}{2}(25-2)$.


\end{exam}

\subsection{A linear time method}

Now we would improve the time complexity of the balanced $2$-maxian problem to linear time.

\begin{lem}\label{endnod}
Let $T$ be a tree. Then an optimal solution of the balanced 2-maxian problem on $T$, exist on the leaf nodes of $T$.
\end{lem}

\begin{prf}
Let $X=\{x_1,x_2\}$ be the solution of balanced 2-maxian problem and $V_1$ and $V_2$ be the sets of vertices that are assigned to $x_1$ and $x_2$, respectively. Let $T_1=(V_1,E_1)$ and $T_2=(V_2,E_2)$ be two subtrees of $T$ which are obtained by deleting edge $e=(v_r,v_s)$ where $v_r\in T_1$ and $v_s \in T_2$. Then $x_1\in T_2$ and $x_2\in T_1$. If either $x_1$ or $x_2$ is not leaf node, then we consider its adjacent. Let $x_1$ be an inner vertex and $u\in T_2$ be its adjacent vertex that is not in the path connecting $x_1$ to $x_2$. In fact, $u$ is a vertex along the path connecting $x_1$ to a leaf nodes. 
Then for all vertices $v_i\in T_1$,
$$d(v_i,u)=d(v_i,x_1)+d(x_1,u)\geq d(u,x_1).$$
Therefore,
$$\sum_{v_i\in T_1}w_id(v_i,u)\geq \sum_{v_i\in T_1}w_id(v_i,x_1).$$
So by relocation of facilities toward leaf nodes the maxian part of the objective function is not decreased while the balancing part remains unchanged.


\end{prf}

As previously mentioned, the optimal solution of $2$-maxian problem is two end vertices of the longest path of the tree. In the following, we show that this property holds for the balanced case.

\begin{lem}\label{notlong}
Let $P$ be the path between two vertices $x_1$ and $x_2$ in the tree $T$.  Let $T_1$ and $T_2$ be two subtrees of $T$ obtained by deleting edge $(v_r,v_s)$ and contain $x_2$ and $x_1$, respectively. If $P$ is not the longest path in the tree $T$, then there exist either a vertex $u\in T_2$ such that $d(u,x_1)\geq d(x_1,x_2)$ or a vertex $u'\in T_1$ such that
$d(u',x_2)\geq d(x_1,x_2).$
\end{lem}

\begin{figure}[h]
  \centering
\includegraphics[width=12cm]{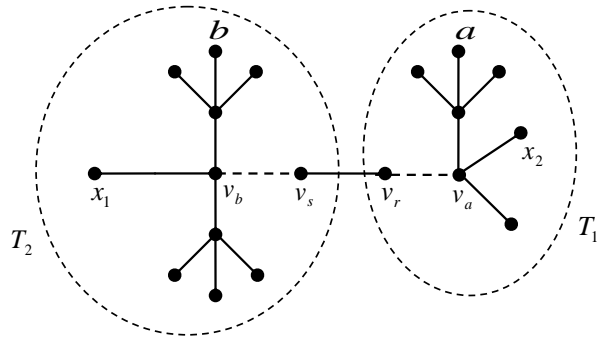}
\vspace{-1.cm}
\caption{Longest and non longest paths on a tree}\label{figlong}
\end{figure}

\begin{prf}
Let $P'$ be the diameter of the tree $T$ and $a$ and $b$ be two end vertices of $P'$. Let $v_a\in P$ and $v_b\in P$ be the closest vertices in $P$ to $a$ and $b$, respectively (see Fig.\ref{figlong}). We consider two cases where the two vertices $a$ and $b$ are in the same or different subtrees. First let $a,b\in T_1$, then either $$d(a,v_a)\geq d(v_a,x_2)~~or~~d(b,v_b)\geq d(v_b,x_2).$$ Otherwise, the path $P'$ is not diameter of the tree. Therefore, either
$$d(a,x_1)=d(a,v_a)+d(v_a,x_1)\geq d(v_a,x_2)+d(v_a,x_1)=d(x_1,x_2),$$
or
$$d(b,x_1)=d(b,v_b)+d(v_b,x_1)\geq d(v_b,x_2)+d(v_b,x_1)=d(x_1,x_2).$$
Now consider the other case where $a\in T_1$ and $b\in T_2$. Then either
$$d(a,v_a)\geq d(v_a,x_2)~~or~~d(b,v_b)\geq d(v_b,x_1).$$
Otherwise, the path $P'$ is not the longest one. So
$$d(a,x_1)\geq d(x_1,x_2),~~or~~d(b,x_2)\geq d(x_1,x_2).$$
The other cases where $a,b\in T_2$ and $a\in T_2,~ b\in T_1$ would be similarly proved.
\end{prf}

\begin{thm}\label{longp}
There is an optimal solution of the balanced 2-maxian problem on the leaf nodes of diameter of $T$.
\end{thm}

\begin{prf}
By Lemma \ref{endnod} there is an optimal solution on the leaf nodes. Let $x_1$ and $x_2$ be two leaf nodes of $T$ that the path connecting them is not diameter of $T$. Let $T_1$ and $T_2$ be two subtrees of $T$ obtained by deleting edge $(v_r,v_s)$ which contain $x_2$ and $x_1$, respectively. Let also the objective function of $x_1$ and $x_2$ respect to this partition be less than or equal to other partitions. The vertices in $T_1$ are allocated to $x_1\in T_2$ and the vertices in $T_2$ are allocated to $x_2\in T_1$ (see Fig. \ref{figlong}). Since the path connecting $x_1$ and $x_2$ is not diameter, by Lemma \ref{notlong} either there is a vertex $u\in T_2$ that $d(u,x_2)\geq d(x_1,x_2)$ or there is a vertex $u'\in T_1$ that $d(u',x_1)\geq d(x_1,x_2).$\\
Without loss of generality, let there exist $u\in T_2$ where $d(u,x_2)\geq d(x_1,x_2)$. Then
$$d(u,v_r)+d(v_r,x_2)=d(u,x_2)\geq d(x_1,x_2)=d(x_2,v_r)+d(v_r,x_1),$$
and consequently
$$d(u,v_r)\geq d(x_1,v_r).$$ Hence, for each $v_i\in T_1$, $$d(v_i,x_1)=d(v_i,v_r)+d(v_r,x_1)\leq d(v_i,v_r)+d(v_r,u)=d(v_i,u).$$
So
$$\sum_{v_i\in T_1}w_id(v_i,u)\geq \sum_{v_i\in T_1}w_id(v_i,x_1).$$
Since the partitions remain unchanged, then the objective function will not be decreased by choosing $u$ instead of $x_1$.
\end{prf}

Note that if $d(v_i,v_j)>0$ for $i,j=1,...,n$ then by Theorem \ref{longp} the optimal solution is two ends of diameter.

Since diameter of a tree can be found in a linear time (see e.g. \cite{H73}), then the optimal solution of balanced 2-maxian problem can be found in $O(n)$ time. However, to calculate the value of objective function we should compute the corresponding objective function by deleting any edge on diameter which would be performed in $O(n^2)$ time. To reduce the time complexity, we first consider the computing objective function on a path.

\begin{lem}\label{cf1f2}
Let $P=v_1,...,v_n$ be a path and $e_1=(v,u)$ and $e_2=(u,s)$ be two adjacent edges on $P$. Let $f_{pmax}^{e_1}(v_1,v_n)$ and $f_{pmax}^{e_2}(v_1,v_n)$ be the objective function values of balanced 2-maxian problem in vertices $v_1$ and $v_n$ by deleting edges $e_1$ and $e_2$, respectively. Also let $Z=\sum_{i=1}^nz_i$ and $v_r$ be the vertex of $P$ so that $\sum_{i=1}^{r-1}z_i<\frac{Z}{2}$ and $\sum_{i=1}^{r}z_i\geq\frac{Z}{2}$. Then
\begin{align} &f_{pmax}^{e_1}(v_1,v_n)-f_{pmax}^{e_2}(v_1,v_n)= \\ &\left\{\begin{array}{ccc}
\lambda w_u(d(u,v_1)-d(u,v_n))+(1-\lambda)2z_u& if &u\in \{v_1,...,v_{r-1}\}\\
\lambda w_u(d(u,v_1)-d(u,v_n))-(1-\lambda)2z_u& if &u\in \{v_{r+1},...,v_{n}\}.\nonumber
       \end{array}\right. \end{align}

\end{lem}
\begin{prf}
For $i=1,2$, let $T_1^{e_i}$ and $T_2^{e_i}$ be the subpaths of $P$ obtained by deleting edge $e_i$. Let  $T_1^{e_i}$ and $T_2^{e_i}$, for $i=1,2$, contain the vertices $v_n$ and $v_1$, respectively. Then $T_1^{e_2}=T_1^{e_1}\setminus\{u\}$ and $T_2^{e_2}=T_2^{e_1}\cup\{u\}$. Therefore,
$$f_{pmax}^{e_1}(v_1,v_n)=\lambda(\sum_{v_i\in T_1^{e_2}}w_id(v_i,v_1)+w_ud(u,v_1)+\sum_{v_i\in T_2^{e_2}}w_id(v_i,v_n)-w_ud(u,v_n))$$
$$-(1-\lambda)|(\sum_{v_i\in T_1^{e_2}}z_i+z_u)-(\sum_{v_i\in T_2^{e_2}}z_i -z_u)|$$
$$=\left\{\begin{array}{ccc}
f_{pmax}^{e_2}(v_1,v_n)+\lambda w_u(d(u,v_1)-d(u,v_n))+2(1-\lambda)z_u& if &u\in \{v_1,...,v_{r-1}\}\\
f_{pmax}^{e_2}(v_1,v_n)+\lambda w_u(d(u,v_1)-d(u,v_n))-2(1-\lambda)z_u& if &u\in \{v_{r+1},...,v_{n}\}.\\
       \end{array}\right.$$
\end{prf}

Let $e_i=(v_i,v_{i+1})$ for $i=1,...,n-1$, then the optimal objective function on the path $P$ can be iteratively computed. First $f_{pmax}^{e_i}(v_1,v_n)$ should be computed for $i=1,r$. Then by using Lemma \ref{cf1f2} the objective function corresponding to deletion other edges on path $P$ will be obtained. Therefore, the total time complexity is $O(n)$.

\begin{figure}\label{figthmneg}
  \centering
\includegraphics[width=12cm]{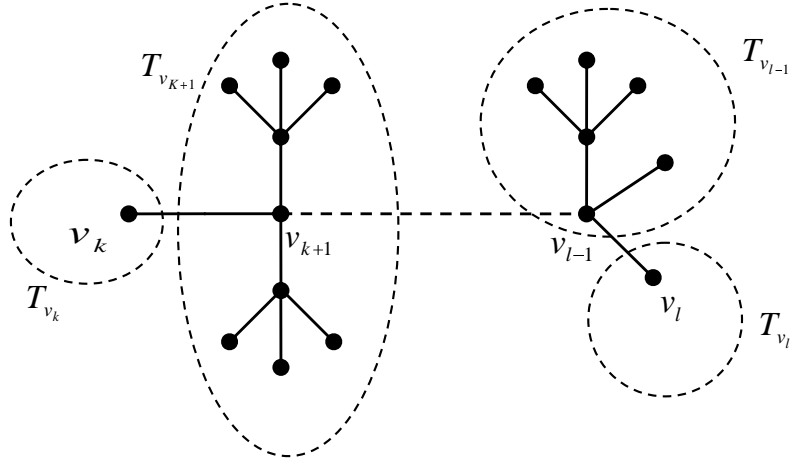}
\vspace{-1.5cm}
\caption{Components of a tree}\label{figthmneg}
\end{figure}

Now consider the tree $T$. Let $P:v_k,...,v_l$ be the diameter of $T$. We create a new path $\widehat{P}$ that the vertices of which are the same as $P$ but the weights and service times of the vertices are varying and defined as follow:
$$\hat{w}_{k}=\sum_{i\in T_{v_{k}}}w_{i},~ \hat{w}_{k+1}=\sum_{i\in T_{v_{k+1}}}w_{i},...,\hat{w}_{l}=\sum_{i\in T_{v_{l}}}w_{i},$$
$$\hat{z}_{k}=\sum_{i\in T_{v_{k}}}z_{i}, ~\hat{z}_{k+1}=\sum_{i\in T_{v_{k+1}}}z_{i},...,\hat{z}_{l}=\sum_{i\in T_{v_{l}}}z_{i}.$$
Where $T_{v_{i}}$, $i=k,k+1,...,l$, are the subtrees of $T$ that obtained by deleting only the edges (not vertices) of $\overline{P}$ from the tree $T$ so that $v_{i}\in T_{v_{i}}$ (see Fig. \ref{figthmneg}). By finding the best deleted edge on $\widehat{P}$, the best deleted edge on the tree $T$ is determined. Therefore, the following theorem is concluded.

\begin{thm}
The balanced 2-maxian problem can be solved in a linear time.
\end{thm}

In the following example the numerical results of the presented methods are given. The results confirm the validity of the findings in the previous sections.

\begin{exam}\label{nofar}

Consider the tree depicted in Fig. \ref{tree17}, where the weights of its vertices are given in the Table \ref{t-2med17}.

\begin{table}[h]
 \centering
\begin{tabular}{c  c  c  c  c  c  c  c  c  c  c  c}
$w_{1}$ & $w_{2}$ &  $w_{3}$ &  $w_{4}$ &  $w_{5}$ &  $w_{6}$ &  $w_{7}$ &  $w_{8}$ &  $w_{9}$ & $w_{10}$ & $w_{11}$ & $w_{12}$ \\
\hline
2 &   4 &   1 &  10 &  3  & 4 &  2 &   2 & 8 &   5 &  5 &   5\\
\hline\hline\\
$w_{13}$ & $w_{14}$ &  $w_{15}$ &  $w_{16}$ &  $w_{17}$ & & && & && \\
\hline
7 & 8 &   3 &  2 & 6 &  &  &  & & & &\\
\hline
\end{tabular}
\caption{The weights of vertices of tree in Fig. \ref{tree17}.}\label{t-2med17}
\end{table}

\begin{figure}[h]
  \centering
\includegraphics[width=12cm]{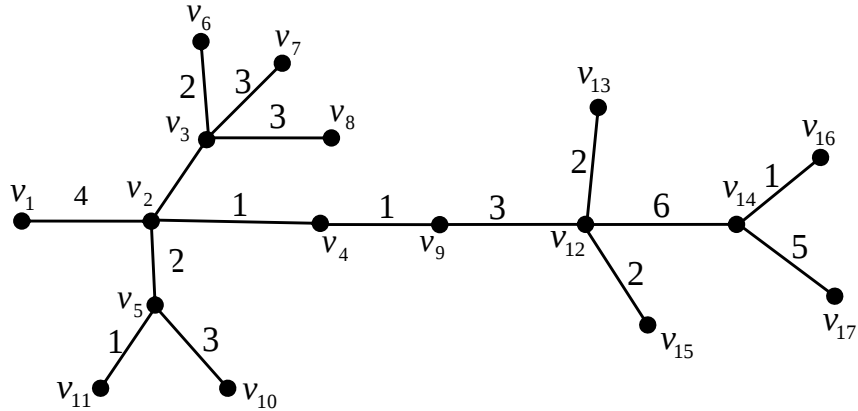}
\caption{A tree with 17 vertices}\label{tree17}
\end{figure}

The solutions of the balanced 2-median problem for varying values of $\lambda$ are presented in Table \ref{t-rmed2-17}. In this table the column with heading $f_5$ indicates the difference between number of clients assigned to each facility. Note that, in the case $\lambda=0.6$, by adding the balanced objective function to median model, we could find a nearly equity solution which its value of the median objective function, i.e. $f_1$, is not considerably increased.

\begin{table}[h]
 \centering
\begin{tabular}{c | c c c c c  c }
\hline
 & $f_1$&$f_5$&deleted edge & $f_{pmed}$ & medians \\
\hline
$\lambda=1$ & 231& 61-16=45&$e=(v_{12},v_{14})$ & 231 &  $\{v_{4},v_{14}\}$\\
$\lambda=0.6$ & 251&46-31=15&$e=(v_{9},v_{12})$ & 156.6 &  $\{v_{2},v_{14}\}$ \\
$\lambda=0.5$&265&39-38=1& $e=(v_{4},v_{9})$ & 133 &  $\{v_{2},v_{12}\}$ \\
$\lambda=0$ &-&39-38=1&$e=(v_{4},v_{9})$ & 1 &  -\\
\hline
\end{tabular}
\caption{The solutions of the balanced 2-median problem for varying values of $\lambda$ on tree in Fig. \ref{tree17}.}\label{t-rmed2-17}
\end{table}

Table \ref{t-rmax2-17} contains the solutions of balanced 2-maxan problem for varying values of $\lambda$. Note that, although for all amounts of $\lambda$ the optimal solution is $\{v_{10},v_{17}$ which are the two end vertices of the diameter. However, the vertices assigned to the facilities are different.  Furthermore, in the cases $\lambda=0.6, 0.5$, by adding the balanced objective function, a nearly equity solution is found, which its value of maxian objective function, i.e. $f_2$, is not considerably decreased.

\begin{table}[h]
 \centering
\begin{tabular}{c | c  c c c c c }
\hline
  & $f_2$&$f_5$& deleted edge & $f_{pmax}$ & optimal solution \\
\hline
$\lambda=1$ &1266&61-16=45& $e=(v_{12},v_{14})$ & 1266 &  $\{v_{10},v_{17}\}$\\
$\lambda=0.6$ &1251&46-31=15& $e=(v_{9},v_{12})$ & 744.6 &  $\{v_{10},v_{17}\}$ \\
$\lambda=0.5$&1251& 46-31=15& $e=(v_{9},v_{12})$ & 618 &  $\{v_{10},v_{17}\}$ \\
$\lambda=0$ &-&39-38=1&$e=(v_{4},v_{9})$ & -1 &  -\\
\hline
\end{tabular}
\caption{The solutions of the balanced 2-maxian problem for varying values of $\lambda$ on tree in Fig. \ref{tree17}.}\label{t-rmax2-17}
\end{table}

\end{exam}

\section{Computational results}
In this section some numerical examples are given for the balanced 2-median and 2-maxian problems. The algorithms were written in MATLAB and tested for 10 randomly generated problems. The arc lengths are generated randomly and taken from the set $[0,5]$.
We assigned the weight 5 to all vertices. We also examined the case that the weights of vertices are generated randomly in the interval $[0,5]$. However, in this case, since for trees with more than 50 nodes, the distances between vertices will be large, then in the most of test problems with varying values of $\lambda$ the functions $f_1$ and $f_2$ dominated on $f_5$. Therefore, the total objective function hasn't considerably changed for $0<\lambda\leq 1$.

The results for varying values of $\lambda$ are presented in Tables \ref{t2medb} and \ref{t2maxb}.
The results show that in some test problems the optimal solutions of classical 2-median and 2-maxian problems are balanced (see Test No. 4). However, in some other test problems, the balancing is depended on $\lambda$. In the most of cases, the difference number of allocated clients to each facility are jumped when $\lambda$ is changed from 0.2 to 0.5. We also examined other values of $\lambda$, but the solutions are not considerably changed.

\begin{table}
  \centering
\begin{scriptsize}
\begin{tabular}{|l|c|c|c|c|c|c|c|c|c|c|}
\hline
&  & $\lambda=0$ & \multicolumn{3}{|c|}{$\lambda=0.2$} & \multicolumn{3}{|c|}{$\lambda=0.5$}& \multicolumn{2}{|c|}{$\lambda=1$}\\ \cline{3-11}
Test\# &n&$f_{pmed}$& $f_1$ &$f_5$ &$f_{pmed}$ & $f_1$ &$f_5$ &$f_{pmed}$&$f_5$ &$f_{pmed}$\\
\hline
1&41&65&1415&65&335&1345&85&715&95&1340\\
2&108&75&8545&175&1849&8565&175&4360&175&8545\\
3&159&165&13580&265&2928&13560&275&6917.5&425&13500\\
4&186&10&14845&10&2977&14845&10&7427.5&10&14845\\
5&243&65&22325&455&4829&22305&475&11390&475&22305\\
6&301&85&36930&595&7862&36930&595&18763&605&36925\\
7&344&190&37500&790&8132&37470&800&19135&800&37470 \\
8&408&310&51695&330&10603&51695&330&26013&330&51695\\
9&463&385&56005&385&11509&25195&1385&27790&1385&54195\\
10&534&500&59450&500&12290&59450&500&29975&500&59450\\
\hline
\end{tabular}
\caption{Results for the balanced 2-median problem. \label{t2medb}}
\end{scriptsize}
\end{table}

\begin{table}
  \centering
\begin{scriptsize}
\begin{tabular}{|l|c|c|c|c|c|c|c|c|c|c|}
\hline
 &  & $\lambda=0$ & \multicolumn{3}{|c|}{$\lambda=0.2$} & \multicolumn{3}{|c|}{$\lambda=0.5$}& \multicolumn{2}{|c|}{$\lambda=1$}\\ \cline{3-11}
Test\#&n&$f_{pmax}$& $f_2$ &$f_5$ &$f_{pmax}$ & $f_2$ &$f_5$ &$f_{pmax}$&$f_5$ &$f_{pmax}$\\
\hline
1&41&-65&5455&65&1039&5505&85&2710&85&5505\\
2&108&-75&28705&75&5681&28705&75&14315&75&28705\\
3&159&-165&47250&265&9238&47260&275&23493&275&47260\\
4&186&-10&46285&10&9249&46285&10&23138&30&46285\\
5&243&-65&83600&445&16364&83625&455&83625&455&83625\\
6&301&-85&112030&85&22338&112030&85&55973&835&112490\\
7&344&-190&113700&190&22588&113700&190&113700&190&113700\\
8&408&-310&177230&320&35190&177260&330&88465&330&177260\\
9&463&-385&157610&385&31214&157610&385&78613&385&157610\\
10&534&-500&172355&500&34071&172355&500&85928&1410&172355 \\
\hline
\end{tabular}
\caption{Results for the balanced 2-maxian problem. \label{t2maxb}}
\end{scriptsize}
\end{table}

\begin{figure}[h]
  \centering
\includegraphics[width=12cm]{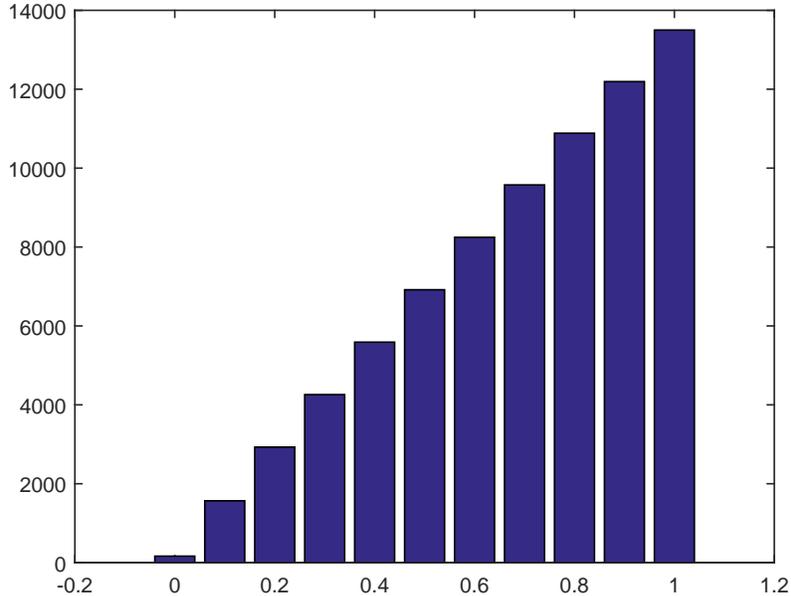}
\caption{The histogram of $f_{pmed}$ respect to $\lambda$ for test NO. 3}\label{histmed3}
\end{figure}

\begin{figure}[h]
  \centering
\includegraphics[width=12cm]{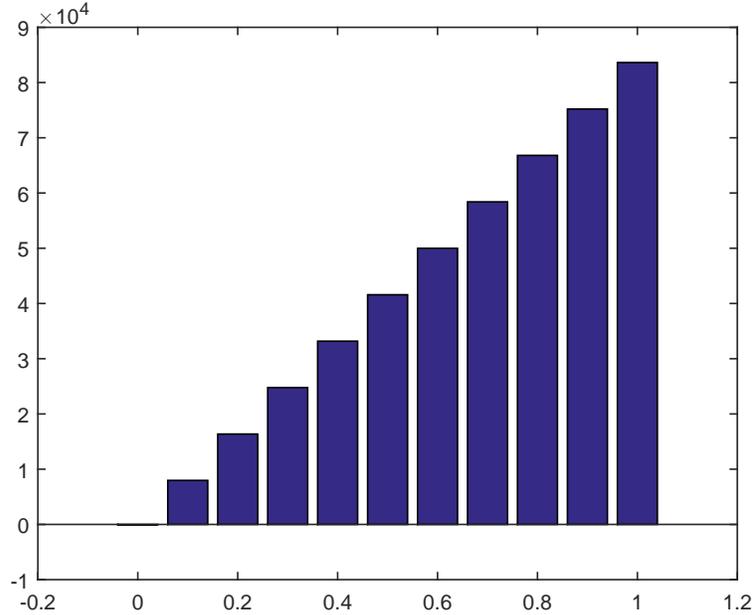}
\caption{The histogram of $f_{pmax}$ respect to $\lambda$ for test NO. 5}\label{histmax5}
\end{figure}

Fig. \ref{histmed3} shows the histogram of the changing values of $f_{pmed}$ for test problem No. 3 respect to $\lambda$. The histogram of the changing values of $f_{pmax}$ for test problem No. 5 respect to $\lambda$, is given in Fig. \ref{histmax5}.

\section{Summary and conclusion}
In this paper, two bi-objective balanced models of the $p$-median and $p$-maxian problems on a tree have been investigated. In the balanced $p$-median problem, the objective function is combination of balance on clients' allocation to the facilities and median problem while in the balanced $p$-maxian problem the objective function is balancing on clients' allocation and maxian problem. Based on edge deletion method an $O(n^2)$ algorithm is presented for the balanced $2$-median problem on a tree. Furthermore, it is shown that the optimal solution of the balanced 2-maxian problem, is the leaf nodes of the diameter of the tree. Then a linear time algorithm is presented to obtain the balanced 2-maxian objective function. To illustrate the algorithms, some numerical examples are given. The results of these examples show that enforcing the balanced objective function to the median and maxian models, causes an almost equitable assignment clients to servers.


\begin{thebibliography}{99}


\bibitem{BP15} Barbati ., Piccolo C., Equality measures properties for location problems, Optimization Letters, 10 (2015), 903-920.

\bibitem{BDTW09}Berman O., Drezner Z., Tamir A., Wesolowsky G.O., Optimal location with equitable loads, Annals of Operations Research, 167 (2009), 307-325.

\bibitem{BFK07} Burkard R.E., Fathali J., Kakhki H.T., The p-maxian problem on a tree, Operations Research Letters, 35 (2007), 331-335.

\bibitem{BCD00} Burkard R.E., \c{C}ela E., Dollani H., 2-Median in trees with pos/neg weights, Discrete Appl. Math., 105 (2000), 51-71.

\bibitem{CG78}Church R.L., Garfinkel R.S., Locating an obnoxious facility on a network, Transportation Science. 12 (1978) 107–118.


\bibitem{EL95}Eiselt H.A., Laporte G. 1995, Objectives in Location Problems. In: Facility Location: A Survey of Applications
and Methods. Ed.: Drezner Z. Springer, Berlin, 151-180.

\bibitem{GH95} Gavalec M., Hudec O., 1995, Balanced location on a graph, Optimization, 35:4, 367-372.

\bibitem{GS95} Gavish B., Sridhar S., Computing the 2-median on tree networks is $O(n\log n)$ time, Networks, 26 (1995), 305-317.

\bibitem{H64}Hakimi S.L., Optimum locations of switching centers and the absolute centers and medians of a graph. Operations Research
12(1964):450-459

\bibitem{H65}Hakimi S.L., Optimum distribution of switching centers in a communication network and some related graph theoretic problems. Operatios Research
13(1965):462–475

\bibitem{H73} Handler G.Y., Minimax location of a facility in an undirected tree networks, Transportation Sci., 7 (1973), 287-293.

\bibitem{HM88} Hansen P. , Moon I.D., Dispersing facilities on a network, Presentation at the TIMS/ORSA Joint National Meeting, Washington D.C.,
April 1988.

\bibitem{KC10} Kang L., Cheng Y., The p-maxian problem on block graphs, Journal of Combinatorial Optimization, 20 (2010), 131-141.

\bibitem{KH79} Kariv O., Hakimi S.~L., An algorithmic approach to network location problems. Part II: p-medians, SIAM J. Appl. Math., 37 (1979), 539-560.

\bibitem{LM14} Landete M., Marin A., Looking for edge-equitable spanning trees, Computers \& Operations Research, 41 (2014), 44-52.

\bibitem{LP13} Lejeune M. A., Prasad S. Y., Effectiveness-equity models for facility location problems on tree networks, networks, 62 (2013), 243-254.

\bibitem{LMP08} Lopez-de-los-Mozosa M.C., Mesaa J.A., Puertoc J., A generalized model of equality measures in network location problems, Computers \& Operations Research, 35 (2008), 651-660.

\bibitem{M11}Marin A., 2011, The discrete facility location problem with balanced allocation of customers. European Journal of Operational Research, 210, 27-38.

\bibitem{MS94}Marsh M.T., Schilling D.A. 1994, Equity Measurement in Facility Location Analysis: A Review and Framework.
European Journal of Operational Research 74 (1): 1-17.

\bibitem{T91}Tamir A., Obnoxious facility location on graphs, SIAM Journal on Discrete Mathematics, 4 (1991) 550–567.

\bibitem{T96} Tamir A., An $O(pn^2)$ algorithm for the p-median and related problems on tree graphs, Operations Research Letters, 19 (1996), 59-64.

\bibitem{T84} Ting S.S., A linear-time algorithm for maxisum facility location on tree networks, Transportation Science, 18 (1984), 76-84.

\bibitem{Z68} Zelinka B., Medians and peripherians of trees, Archivum Mathematicum,  4 (1968), 87-95.

\end{thebibliography}
\end{document}